\documentclass[oribibl,envcountsame,runningheads]{llncs}

\usepackage{amsmath}
\usepackage{amssymb}
\usepackage{latexsym}

\let\url\relax
\usepackage[colorlinks=false,bookmarks=true,pdftitle={Quantified Propositional Gödel Logic},pdfauthor={Matthias Baaz, Agata Ciabattoni, and Richard Zach},pdfkeywords={many-valued logic, intuitionistic fuzzy logic, proof theory, Gödel-Dummett logic},dvips]{hyperref}

\spnewtheorem{defn}[theorem]{Definition}{\bfseries}{\rmfamily}

\let\impl\supset
\let\nmodels\nvDash

\def\qd{\mathop\mathrm{qd}}
\def\nfp{^{(p)}}
\def\log#1{\ensuremath{\mathbf{\mathsf{#1}}}}
\def\LC{\ensuremath{\mathbf{LC}}}
\def\Var{\mathop\mathit{Var}}
\def\OP{\mathrm\mathit{OP}}
\def\qf{\mathrm{qf}}
\def\nf{\mathrm{nf}}
\def\qp{\mathrm{qp}}
\let\up\uparrow
\let\down\downarrow

\newcommand{\Q}{\mathsf{Q}}

\def\oper{\hbox{\Large$\circ$}}

\newcommand{\ifff}{\leftrightarrow}
\newcommand{\prov}{\QG \vdash}
\newcommand{\QG}{\ensuremath{\log{QG}^\qp_\uparrow}}
\newcommand{\G}{\ensuremath{\mathbf{G}^\qp_\uparrow}}

\title{Quantified Propositional G\"odel
Logics\protect\footnotetext{\protect\href{http://e-math.ams.org/msc/}
{\emph{2000 Mathematics Subject Classification:}} 
Primary 03B50; Secondary 03B55.}\thanks{Research supported by the Austrian Science Fund under
grant P--12652~MAT}}

\author{Matthias Baaz\inst{1} \and 
Agata Ciabattoni\inst{1}$^,$\thanks{Research supported by EC Marie Curie
fellowship HPMF--CT--19 99--00301} \and
\protect\href{http://www.logic.at/people/zach/}{Richard Zach}\inst{2}}

\authorrunning{Matthias Baaz, Agata Ciabattoni, and Richard Zach}

\institute{Institut f\"ur Algebra und Computermathematik E118.2,\\
Technische Universit\"at Wien,
A--1040 Vienna, Austria,
\email{[\protect\href{mailto:baaz@logic.at}{baaz}, 
\protect\href{mailto:agata@logic.at}{agata}]@logic.at}
\and
Institut f\"ur Computersprachen E185.2,\\
Technische Universit\"at Wien,
A--1040 Vienna, Austria,
\email{\protect\href{mailto:zach@logic.at}{zach@logic.at}}
}

\makeatletter\newcommand{\ps@ref}{\addtolength{\headheight}{3ex}
\addtolength{\topmargin}{-5ex}
\addtolength{\headsep}{2ex}
\renewcommand{\@evenhead}{}
\renewcommand{\@evenfoot}{\footnotesize\renewcommand{\arraystretch}{.8}{
\noindent\begin{tabular}{@{}l} Parigot, Michel and Andrei Voronkov
(eds.), \emph{Logic for Programming and}\\ \emph{Automated Reasoning.
7th International Conference, LPAR 2000.}\\
Proceedings. \href{http://www.springer.de/cgi-bin/search_book.pl?isbn=3-540-41285-9}{LNCS 1955}. Springer,
Berlin, 2000. pp.~240--256\\ \copyright Springer-Verlag Berlin Heidelberg New York
2000
\end{tabular}
}\hfill} } 
\def\ps@headings{\let\@mkboth\@gobbletwo
   \let\@oddfoot\@empty\let\@evenfoot\@empty
   \def\@oddhead{\normalfont\small\rlap{\thepage}\hspace{\headlineindent}%
                  \leftmark\hfil}
   \def\@evenhead{\normalfont\small\hfil\rightmark\hspace{\headlineindent}%
                 \llap{\thepage}}
   \def\chaptermark##1{}%
   \def\sectionmark##1{}%
   \def\subsectionmark##1{}}

\makeatother

\begin{document} 
\bibliographystyle{abbrv}
\maketitle\thispagestyle{ref}\pagestyle{headings}


\begin{abstract}
It is shown that \G, the quantified propositional G\"odel logic based
on the truth-value  set $V_\uparrow = \{1 - 1/n : n \ge 1\}\cup\{1\}$,
is decidable.  This result is obtained by reduction to B\"uchi's
theory S1S.  An alternative proof based on elimination of quantifiers
is also given, which yields both an axiomatization and a
characterization of \G{} as the intersection of all finite-valued
quantified propositional G\"odel logics.
\end{abstract}


\section{Introduction}

In 1932, G\"odel~\cite{goedel} introduced a family of finite-valued
propositional logics to show that intuitionistic logic does not have a
characteristic finite matrix.  Dummett~\cite{dummett} later
generalized these to an infinite set of truth-values, and showed that
the set of its tautologies {\bf LC} is axiomatized by intuitionistic
logic extended by the linearity axiom $(A \impl B) \lor (B \impl A)$.
G\"odel-Dummett logic naturally turns up in a number of different
areas of logic and computer science.  For instance, Dunn and
Meyer~\cite{DM} pointed out its relation to relevance logic;
Visser~\cite{visser} employed it in investigations of the provability
logic of Heyting arithmetic; Pearce used it to analyze inference in
extended logic programming \cite{PearceD:99}; and eventually it was
recognized as one of the most important formalizations of fuzzy logic
\cite{hajek}.

The propositional G\"odel logics are well understood: Any infinite set
of truth-values characterizes the same set of tautologies.  {\bf
LC}~is also characterized as the intersection of the sets of
tautologies of all finite-valued G\"odel logics
$\mathbf{G}_k$~\cite{dummett}, and as the logic determined either by
linearly ordered Kripke frames or linearly ordered Heyting
algebras~\cite{horn}.

When G\"odel logic is extended beyond pure propositional logic,
however, the situation is more complex.  For
the cases of propositional entailment and extension to first-order
validity, infinite truth-value sets with different order types
determine different logics with different properties.  There are
infinitely many sets of truth values which give rise to distinct
logics.  As an example, consider the truth-value sets
\begin{eqnarray*}
V_\infty & = &[0, 1]\\
V_\downarrow & =& \{0\} \cup \{1/n : n \ge 1\} \\
V_\uparrow & = &\{1\} \cup \{1-1/n : n \ge 1\}\\
V_k & = & \{1\} \cup \{1-1/n : n = 1, \dots, k-1\}
\end{eqnarray*}
Propositional entailment with respect to $V_\infty$ is compact, but
not with respect to $V_\downarrow$ or $V_\up$.  If a formula $A$ is
entailed by a set $\Gamma$ with respect to $V_k$ for every $k$, then
it is also entailed with respect to $V_\up$, but not necessarily with
respect to $V_\infty$ or $V_\down$ \cite{BaazZach}.  Similarly, the
first-order logic based on $V_\infty$ is axiomatizable (this is
Takeuti and Titani's intuitionistic fuzzy logic \cite{takeuti}), while
those based on $V_\up$ and $V_\downarrow$ are not
\cite{baaz-leitsch-zach}.  The first-order G\"odel logic based on
$V_\up$ is the intersection of all finite-valued
first-order G\"odel logics.

Another interesting generalization of propositional logic is obtained
by adding quantifiers over propositional variables.  In classical
logic, propositional quantification does not increase expressive power
per se.  It does, however, allow expressing complicated properties
more naturally and succinctly, e.g., satisfiability and validity of
formulas are easily expressible within the logic once such quantifiers
are available.  This fact can be used to provide efficient proof
search methods for several non-monotonic reasoning
formalisms~\cite{EETW}.  

For G\"odel logic the increase in expressive power is witnessed by the
fact that statements about the topological structure of the set of
truth-values (taken as infinite subsets of the real interval $[0,1]$)
can be expressed using propositional
quantifiers~\cite{baaz-veith-98b}.  In~\cite{baaz-veith-98b} it is
also shown that there is an uncountable number of different quantified
propositional infinite-valued G\"odel logics.  The same paper
investigates the quantified propositional G\"odel logic
$\mathbf{G}^\qp_\infty$ based on the set of truth-values $[0, 1]$,
which was shown to be decidable. It is of some interest to
characterize the intersection of all finite-valued quantified
propositional G\"odel logics.  As was pointed out in
\cite{baaz-veith-98b}, $\mathbf{G}^\qp_\infty$ does not provide such a
characterization.

In this paper we study the quantified propositional G\"odel logic $\G$
based on the truth-value set~$V_\uparrow$.  We show that \G{} is
decidable.  In general, it is not obvious that a quantified
propositional logic is decidable or even axiomatizable. For instance,
neither the closely related quantified propositional intuitionistic
logic, nor the set of valid first-order formulas on the truth-value
set $V_\uparrow$ are r.e.  Although our result can be obtained by
reduction to B\"uchi's monadic second order theory of one successor
S1S \cite{buechi}, we also give a more informative proof based on
elimination of propositional quantifiers.  This proof allows us to
characterize \G{} as the intersection of all finite-valued quantified
propositional G\"odel logics, and moreover yields an axiomatization
of~\G.

A remark is in order about the relationship between the approach taken
here using truth-value semantics and Kripke semantics. As was pointed
out above, {\bf LC} is often defined as the propositional logic of
linearly ordered Kripke frames.  In Kripke semantics, quantified
propositional \LC{} would then result by adding quantifiers over
propositions (subsets of the set of worlds closed under
accessibility).  Here different classes of linear Kripke structures
which all define {\bf LC} in the pure propositional case in general do
not define the same quantified propositional logic.  In particular,
the logic obtained by just taking Kripke models of order type $\omega$
is not the same as that defined by the class of all finite linear
orders.  It follows from the results of this paper that the logic of
all finite linear Kripke structures coincides with~\G.


\section{G\"odel Logics} \label{sec:goedel}

\paragraph{Syntax.}
We work in the language of propositional logic containing a countably
infinite set $\Var = \{p, q, \ldots \}$ of (propositional)
variables, the constants $\bot, \top$, as well as the connectives
$\land, \lor$, and $\impl$.  Propositional variables and constants are
considered atomic formulas.  Uppercase letters will serve as
meta-variables for formulas.  If $A(p)$ is a formula containing the
variable $p$ free, then $A(X)$ denotes the formula with all
occurrences of the variable $p$ replaced by the formula $X$. $Var(A)$
is the set of variables occurring in the formula $A$.
We use the abbreviations $\neg A$ for $A \supset \bot$
and  $A \ifff B$ for $(A \impl B) \land (B \impl A)$.

\paragraph{Semantics.}
The most important form of G\"odel logic is defined over the real unit
interval $V_\infty = [0,1]$; in a more general framework, the
truth-values are taken from a set $V$ such that $\{0,1\} \subseteq V
\subseteq [0,1]$.  In the case of $k$-valued G\"odel logic ${\bf
G}_k$, we take $V_k = \{1 - 1/i : i = 1, \ldots, k-1\} \cup \{1 \}$.
The logic we will be most interested in is based on the set $V_\up =
\{1 -1/i : i \ge 1\} \cup \{1\}$.

A \emph{valuation}~$v\colon \Var \to V$ is an assignment of values in
$V$ to the propositional variables. It can be extended to formulas using
the following truth functions introduced by G\"odel~\cite{goedel}:
\[
\begin{array}{cc}
\begin{array}{rcl}
v(\bot)                         &=& 0 \\
v(\top)                         &=& 1 \\
v(A \land B)     &=& \min(v(A), v(B))
\end{array} &
\begin{array}{rcl}
v(A \lor B)      &=& \max(v(A), v(B))\\
v(A \supset B)  &=&
 \begin{cases}
   1            & {\rm if\ } v(A) \leq v(B) \\
   v(B)         & {\rm otherwise}
  \end{cases}
\end{array}
\end{array}
\]
A formula $A$ is a \emph{tautology} over a truth-value set $V
\subseteq [0,1]$ if for all valuations $v\colon \Var \to V$, $v(A) =
1$.  The \emph{propositional logics} $\mathbf{LC}$, $\mathbf{G}_\up$
and $\mathbf{G}_k$ are the sets of tautologies over the corresponding
truth value sets, e.g., $\mathbf{LC} = \mathbf{G}_\infty = \{A : A
\textrm{ a tautology over }V_\infty\}$.  We also write $\mathbf{G}
\models A$ for $A \in \mathbf{G}$ ($\mathbf{G} \in \{\mathbf{LC},
\mathbf{G}_\up, \mathbf{G}_k\}$).

It is easily seen that $\mathbf{LC} \supseteq \mathbf{G}_\up \supseteq
\mathbf{G}_k$.  Dummett \cite{dummett} showed that $\mathbf{LC} =
\mathbf{G}_\up$ and that $\mathbf{LC} = \bigcap_{k\ge 2}
\mathbf{G}_k$.

The abbreviation $A \prec B$ for $(A \impl B) \land ((B \impl A) \impl
A)$ will be used extensively below.  It expresses strict linear order
in the sense that
\[
 v(A \prec B) = 
\begin{cases} 
1 &  \mathrm{if\ } v(A) < v(B) \mathrm{\ or\ } v(B) = 1\\
\min(v(A),v(B)) & {\rm otherwise} \end{cases}
\]

\paragraph{Propositional Quantification.}

In {\em classical} propositional logic we define $(\exists p) A(p)$ by
$A(\bot) \lor A(\top)$ and $(\forall p) A(p)$ by $A(\bot) \land
A(\top)$.  In other words, propositional quantification is
semantically defined by the supremum and infimum, respectively, of
truth functions (with respect to the usual ordering ``$0 < 1$'' over
the classical truth-values $\{0,1\}$).  This can be extended to
G\"odel logic by using {\em fuzzy quantifiers}. Syntactically, this
means that we allow formulas $(\forall p) A$ and $(\exists p) A$ in
the language.  Free and bound occurrences of variables are defined in
the usual way.  Given a valuation $v$ and $w \in V$, define $v[w/p]$
by $v[w/p](p) = w$ and $v[w/p](q) = v(q)$ for $q \not\equiv p$.  The
semantics of fuzzy quantifiers is then defined as follows:
\[
v((\exists p) A) = \sup \{ v[w/p](A) : w \in V \}  \hspace{4ex}
v((\forall p) A) = \inf \{ v[w/p](A) : w \in V\}
\]
When we consider quantifiers, $V$ has to be closed under 
infima and suprema, since otherwise truth values for quantified
formulas are not defined. 

We also add the additional unary connective $\oper$ to the
language. The truth function for $\oper$ is given by $v(\oper A) =
v((\forall p) ((p \impl A) \vee p))$.  In \G, this makes
\[
v(\oper A) =\begin{cases} 1 & {\rm if\ } v(A) = 1 \\
                  1 - \frac{1}{n+1} & {\rm if\ } v(A) = 1 - \frac{1}{n} \cr
               \end{cases}
\]
We abbreviate $\oper \dots \oper A$ ($n$ occurrences of $\oper$) by
$\oper^n A$.

Using the above definitions, it is straightforward to extend the
notion of tautologyhood to the new language.  We write $\G$
($\mathbf{G}_\infty^\qp$, $\mathbf{G}_k^\qp$) for the set of
tautologies in the extended language over $V_\up$ ($V_\infty$, $V_k$).

We will show below that every quantified propositional formula is
equivalent in \G{} to a quantifier-free formula, which in general can
contain~$\oper$.  $\oper A$ itself (or the equivalent formula
$(\forall p)((p \impl A) \lor p)$), however, is not in general
equivalent to a quantifier-free formula not containing~$\oper$.
Inspection of the truth tables shows that a quantifier-free formula
containing only the variable~$q$ takes one of $0$, $v(q)$, or 1 as its
value under a given valuation $v$, and thus no such formula can
define~$\oper q$.


\section{Hilbert-style Calculi}

All the calculi we consider are based on the following set of axioms:
\[
\begin{array}{l@{\quad}l@{\qquad}l@{\quad}l}
\mathrm{I1} & A \impl (B \impl A) &
     \mathrm{I7} & (A \land \neg A) \impl B\\
\mathrm{I2} & (A \land B) \impl A &
    \mathrm{I8} & (A \impl \neg A) \impl \neg A\\
\mathrm{I3} & (A \land B) \impl B &
    \mathrm{I9} & \bot \impl A\\
\mathrm{I4} & A \impl (B \impl (A \land B)) &
    \mathrm{I10} & A \impl \top\\
\mathrm{I5} & A \impl (A \lor B) &
\mathrm{I11} & (A \impl (B \impl C)) \impl ((A \impl B) \impl (A \impl C))\\
\mathrm{I6} & B \impl (A \lor B) & 
\mathrm{I12} & ((A \impl C)\land (B \impl C)) \impl ((A \lor B) \impl C)\\
\end{array}
\]
These axioms, together with the rule of modus ponens, define the
system \log{IPC} that is sound and complete for intuitionistic
propositional logic.  The system~\log{LC} is obtained by adding
to~\log{IPC} the linearity axiom 
\[
\mathrm{LC} \quad (A \impl B) \lor
(B \impl A).
\] 
It is well known \cite{dummett} that \log{IPC} and
\log{LC} are sound for all propositional G\"odel logics, and that
\log{LC} is complete for all infinite-valued propositional G\"odel
logics.  We will make frequent use of this fact below, and omit
derivations of formulas which are (instances of) quantifier- and
$\oper$-free tautologies in $\mathbf{G}_\up$.  These omissions are
indicated by pointing out that the formula follows already in \log{LC}
or \log{IPC}.  In particular, familiar inference patterns such as the
chain rule or case distinction are derivable in \log{LC} and its
extensions.

When we turn to quantified propositional logics, a natural
system~$\log{IPC}^\qp$ to start with is obtained by adding to~\log{IPC}
the following two axioms:
\[
{\impl}{\exists} \quad A(C) \impl (\exists p) A(p) \qquad\qquad
{\impl}{\forall} \quad (\forall p) A(p) \impl A(C) 
\]
and the rules:
\[
\frac{A(p) \impl B\nfp}{(\exists p)A(p) \impl B\nfp}
\mathrm{R}{\exists}
\quad \quad \quad
\frac{B\nfp \impl A(p)}{B\nfp \impl (\forall p) A(p)}
\mathrm{R}{\forall}
\]
where for any formula $C$, the notation $C\nfp$ indicates that $p$
does not occur free in $C$, i.e., $p$ is a (propositional) {\em
eigenvariable}. 

Let $\QG$ be the system obtained by adding to~$\log{IPC}^\qp$ the
axioms (LC),
\[
{\forall}{\lor}\qquad (\forall p) [A \lor B(p))] \impl [A \lor (\forall p) B(p)]
\]
where $p \notin A$, and the following:
\[
\begin{array}{l@{\quad}l@{\quad}l@{\quad}l}
\mathrm{G1} & \oper (A \impl B) \ifff (\oper A \impl \oper B) &
\mathrm{G4} & (A \impl \oper B) \impl ((A \impl C) \vee (C \impl B)) \\
\mathrm{G2} & A \prec \oper A &
\mathrm{G5} & (A \ifff \bot) \vee (\exists p) (A \ifff \oper p) \\
\mathrm{G3} & (\oper A \impl \oper B) \impl ((A \impl B) \vee \oper B) &
\mathrm{G6} & (A \prec B) \impl (\oper A \impl B)
\end{array}
\]
\begin{proposition}
The system $\QG$ is sound for $\mathbf{G}_k^\qp$ and~$\G$.
\end{proposition}

\begin{proof}  It is easily seen that the rules of inference preserve 
validity. For instance, if $B \impl A(p)$ is valid, then, for any
valuation~$v$, $v[w/p](B) \le v[w/p](A(p))$ where $w \in V$.  If $p$
does not occur in $B$, then $v(B) = v[w/p](B)$ and we have $v(B) \le
\inf \{v[w/p](A(p)) : w \in V\}$.  That $\log{LC}$ is sound for
arbitrary G\"odel logics was shown in~\cite{dummett}.  The tedious but
straightforward verification that the remaining axioms
(${\lor}{\forall}$) and (G1)--(G6) are valid is left to the reader.
\end{proof}

\begin{remark}
In \cite{baaz-veith-98b} it was shown that a system sound and complete
for $\mathbf{G}_\infty^\qp$, the quantified propositional G\"odel
logic based on the truth-value set $[0, 1]$, is obtained by
extending~$\log{IPC}^\qp$ with (LC), (${\lor}{\forall}$) and the axiom
\[
(\forall p)[(A\nfp \impl p) \lor (p \impl B\nfp)]
                        \impl \ (A\nfp \impl B\nfp).
\]
This schema is not valid in \G{} (it comes out $=0$ under any $v$ with
$v(A) = 1/2$ and $v(B) = 0$).  On the other hand, it is easy to see
that $v(\oper A) = v(A)$ in $V_\infty$, and hence axiom (G2) is not
valid in $\mathbf{G}_\infty^\qp$.  Thus neither of
$\mathbf{G}_\infty^\qp$ and \G{} is included in the other.  This is in
contrast to the situation in propositional entailment and first-order
logic, where $V_\infty$ defines the smallest G\"odel logic and is
included in all others.
\end{remark}


\section{Decidability}

In this section we prove that $\G$ is decidable. This is done by
defining a reduction of tautologyhood in \G{} to S1S, the monadic
theory of one successor, which was shown to be decidable
by~B\"uchi~\cite{buechi}.

S1S is the set of second-order formulas in the language with
second-order quantification restricted to monadic set variables $X$,
$Y$, \dots{} with one unary function $'$ (successor) which are true in
the model~$\langle \omega, '\rangle$.  For the purposes of this
section we consider $\oper A$ to be an abbreviation of $(\forall p)
((p \impl A) \vee p)$.

Suppose $A$ is a quantified propositional formula, and
$B$ is a formula in the language of S1S with only $x$ free.
Let $TV(B(x))$ abbreviate $(\forall z)(B(z') \impl B(z))$.
We define $A^x$ by:
\begin{eqnarray*}
p^x & = & X_p(x) \\
\bot^x & = & X_\bot(x) \\
\top^x & = & (\forall z)(z=z) \\
(B \land C)^x & = & B^x \land C^x \\
(B \lor C)^x & = & B^x \lor C^x\\
(B \impl C)^x & = & (\forall y)(B^y \impl C^y) \lor
 (\exists y)(B^y \land \neg C^y) \land C^x \\
(\forall p)B^x & = & (\forall X_p)(TV(X_p(x)) \impl B^x) \\
(\exists p)B^x & = & (\exists X_p)(TV(X_p(x)) \land B^x)
\end{eqnarray*}
Consider the following reduction:
\[
\Phi(A)  = (\forall X_\bot)((\forall x)\neg X_\bot(x) \impl (\forall x)A^x)
\]
The idea behind this is to correlate truth-values in $V_\up$ with
subsets of~$\omega$ which are closed under predecessor, i.e.,
predicates in
\[
TV = \{P \subseteq \omega : \textrm{if\ } n \in P \textrm{\ then\ } m \in P \textrm{\ for all\ } m \le n\}.
\]
Under this correlation, $1$ corresponds to $\omega$, and $1 - 1/n$
corresponds to $\{1, \ldots, n\}$.  

Let $s$ be an interpretation of the language of S1S, mapping variables
to elements or subsets of~$\omega$.  We denote by $s[n/x]$ the
interpretation which is just like $s$ except that it assigns $n$ to
$x$.  Then $TV(A(x))$ obviously expresses the condition that the
predicate~$A(x)[s] = \{n : S1S \models A(x) [s[n/x]]\}$ defined by
$A(x)$ in $s$ is closed under predecessor.
If a monadic predicate~$P$ is closed under predecessor, we define its
truth value by
\[ 
tv(P) = \sup \{1 - \frac{1}{n} : 1^n \in P\}.
\]
Conversely, every truth-value~$v \in V_\uparrow$ corresponds to a
monadic predicate
\[
mp(v) = \begin{cases}\{k : k \le n\} & \textrm{if\ } v = 1 - 1/n\\ 
               \omega & \textrm{if\ } v = 1.
\end{cases}
\]
Note that for $P, Q \in TV$, $P \subseteq Q$ iff $tv(P)
\le tv(Q)$, and conversely, for $v, w \in V_\uparrow$, $v \le w$ iff
$mp(v) \subseteq mp(w)$.  

\begin{lemma}\label{gifs1s}
Let $v$ be a valuation and $s$ be the interpretation defined by
$s(X_p) = mp(v(p))$ and $s(X_\bot) = \emptyset$.  Then we have
$tv(A^x[s]) = v(A)$.
\end{lemma}

\begin{proof}
By induction on the complexity of $A$.  The claim is obvious for atomic 
formulas, conjunction and disjunction. If $A \equiv B \impl C$ we have
to distinguish two cases. Suppose first that $v(B) \le v(C)$. 
By induction hypothesis,
$B^x[s] = mp(v(B)) \subseteq mp(v(C)) = C^x[s]$, and hence the first
disjunct in the definition of $(B \impl C)^x$ is true.  Thus $(B \impl
C)^x$ defines $\omega$ and $tv((B\impl C)^x[s]) = 1$.  Now
suppose that $v(B) > v(C)$.  Then $tv(B^x[s]) \supsetneq tv(C^x[s])$,
$S1S \nmodels (\forall y)(B^y \impl C^y)\ [s]$ and $S1S \models
(\exists y)(B^y \land \neg C^y)\ [s]$, and thus $(B\impl C)^x[s] =
C^x[s]$.

If $A \equiv (\exists p)B$, let $v[w/p]$ be the
valuation which is just like $v$ except that $v[w/p](p) = w$, and let
$s[mp(w)/X_p]$ be the corresponding interpretation which is like $s$
except that it assigns $mp(w)$ to $X_p$.  

By induction hypothesis, $tv(B^x[s[mp(w)/X_p]]) = v[w/p](B)$.  We
again have two cases.  Suppose first that $\sup \{v[w/p](B) : w \in
V_\uparrow\} = 1 - 1/n$.  For all $m > n$, $S1S \nmodels B^x [m/x,
mp(w)/X_p]$, since $v[w/p](B^x) < 1-1/m$ by induction hypothesis.  On
the other hand, $S1S \models TV(P_p) \impl B^x\ [s[k/x,
mp(1-1/n)/P_p]]$ for all $k\le n$, and so $tv((\exists p)B^x[s]) = 1 -
1/n$.  Now consider the case where $\sup \{v[w/p](B) : w \in
V_\uparrow\} = 1$.  Here there is no bound $n$ on the 
the members of sets defined by $B^x[s[mp(w)/X_p]]$ where $w \in
V_\uparrow$.  Hence, $mp((\exists p)B)^x[s]) = \omega$ and
$tv((\exists p)B^x[s]) = 1$.

The case $A \equiv (\forall p)B$ is similar.\qed
\end{proof}

\begin{lemma}\label{s1sifg}
Let $s$ be an interpretation with $s(X_\bot) = \emptyset$ and $s(X_p)
\in TV$. Let $v$ be defined by $v(p) = tv(s(X_p))$. Then $A^x[s] \in
TV$, and $v(A) = tv(A^x[s])$.
\end{lemma}

\begin{proof}
By induction on the complexity of $A$.  The claim is again trivial for
atomic formulas, conjunctions or disjunctions.  If $A \equiv B
\impl C$, two cases occur.  If $S1S \models (\forall y)(B^y \impl
C^y)$, then $B^y[s] \subseteq C^y[s]$.  By induction hypothesis, $v(B)
\le v(C)$, and hence $v(B \impl C) = 1 = tv((B\impl C)^x[s])$.
Otherwise, for some $n$ we have $n \in B^y[s]$ but $n \notin
C^y[s]$. So $(\exists y)(B^y \land \neg C^y)$ must be true and the
predicate defined is the same as~$C^y[s]$.

Now for the case $A \equiv (\exists p)B$: If $S1S \models (\exists
X_p)(TV(X_p) \impl B^x)[s[n/x]]$, then there is a prefix closed
witness $P$ so that $S1S \models B^x[s[n/x,P/X_p]]$.  By induction
hypothesis, $B^x[s[P/X_p]] \in TV$, and hence $S1S \models
TV(X_p) \impl B^x\ [s[m/x,P/X_p]]$ for all $m \le n$, and thus
$((\exists p)B)^x[s] \in TV$ as well.

Consider $N = ((\exists p)B)^x[s]$.  First, suppose that $\sup N = k$.
That means that for some $P \in TV$, $1^k \in B^x[s[P/X_p]]$, and for
no $Q \in TV$ and no $j > k$, $j \in B^x[s[Q/X_p]]$.  By induction
hypothesis, $v[tv(P)/p](B) = 1 - 1/k$ and for all $w \in V_\uparrow$,
$v[w/p](B) \le 1 - 1/k$.  Hence $v((\exists p)B) = 1 - 1/k$.

If $\sup N$ does not exist, for each $k$ there is a witness $Q_k \in
TV$ with $k \in B^x[s[Q_k/X_p]]$.  By induction hypothesis, for each
$k$ we have $v[tv(Q_k)/p](B) \ge 1 - 1/k$, and so $v((\exists p)B) =
1$.

The case $A \equiv (\forall p)B$ is similar.\qed
\end{proof}

\begin{theorem} \G{} is decidable. \end{theorem}
\begin{proof}
If there is a valuation $v$ such that $v(A)
< 1$, then by Lemma~\ref{gifs1s} there is an $s$ with $s(P_\bot) =
\emptyset$ and $n$ so that $n \notin A^x[s]$, and
hence $S1S \nmodels \Phi(A)$.

Conversely, suppose $S1S \nmodels \Phi(A)$.  We may assume, without
loss of generality, that all propositional variables in $A$ are bound.
Then there is an interpretation $s$ with $X_\bot(x)[s] = \emptyset$ so
that some $n \notin A^x[s]$.  By Lemma~\ref{s1sifg}, $A^x[s] \in TV$.
Hence, if $n \notin A^x[s]$, then $k \notin A^x[s]$ for all $k\ge n$,
and, also by Lemma~\ref{s1sifg}, $v(A) = tv(A^x[s]) < 1$.

Thus a formula $A$ is a tautology in \G{} iff $S1S \models \Phi(A)$.
The claim follows by the decidability of $S1S$.\qed
\end{proof}

\section{Properties and Normal Forms}

In this section we introduce suitable normal forms for formulas of
$\QG$ and prove some useful properties of $\QG$. These results will
be crucial in the proof of the elimination of quantifiers.  

\begin{proposition}
\label{basicpr}
\begin{enumerate}
\item $\prov (A \impl B) \impl (\oper A \impl \oper B)$
\item $\prov \oper (A \wedge B) \ifff (\oper A \wedge \oper B)$
\item $\prov \oper (A \vee B) \ifff (\oper A \vee \oper B)$
\end{enumerate}
\end{proposition}

\begin{proof}
(1) From (G2) we have $(A \impl B) \impl \oper(A \impl B)$, which,
    together with the left-to-right direction of (G1) yields the
    result.

(2) The left-to-right implication immediately follows from 
axioms (I2) and (I3) together with Prop.~\ref{basicpr}(1).   
For the converse, replace $B$ by $B \impl (A \land B)$
in Prop.~\ref{basicpr}(1) and use (I4) to derive 
$\oper A \impl \oper (B \impl (A \wedge B))$. 
Then, using (G1), one has  $\oper A \impl (\oper B \impl 
\oper (A \wedge B))$. The claim follows by \log{IPC}.

(3) In $\log{LC}$, we have $(A \lor B) \ifff (A \impl B) \impl B)
\land (B \impl A) \impl A)$.  Replacing $A$ by $\oper A$ and $B$ by
$\oper B$, we have $(\oper A \lor \oper B) \ifff (\oper A \impl \oper
B) \impl \oper B) \land (\oper B \impl \oper A) \impl \oper A)$. The
result follows using (G1) and \log{IPC}.\qed
\end{proof}

\begin{proposition}\label{equivsubst}
\begin{enumerate}
\item If $p$ does not occur boind in $C(p)$, then
\[ \prov (\forall \bar q)(A \ifff B) \impl (C(A) \impl C(B)) \]
where $\bar q$ are the propositional variables occurring free in $A$
and $B$.

\item If $C(p)$ is quantifier-free, we also have
\[ \prov (A \ifff B) \impl (C(A) \impl C(B)) \]
\end{enumerate}
\end{proposition}

\begin{proof}
By induction on the complexity of $C$. Cases for $\land$, $\lor$, and
$\impl$ are easy.  If $C(p) \equiv \oper D(p)$, we use the induction
hypothesis and Prop.~\ref{basicpr}(1).  If $C(p) \equiv (\exists
r)D(p, r)$, we argue:
\begin{eqnarray*}
(1) & (\forall \bar q)(A \ifff B) \impl (D(A, r) \impl D(B, r)) & \textrm{by IH}\\
(2) & \qquad ((\forall \bar q)(A \ifff B) \land D(A, r)) \impl D(B, r)) & \textrm{(1), \log{IPC}}\\
(3) & D(B, r) \impl (\exists r)D(B, r) & {\impl}{\exists}\\
(4) & (\forall \bar q)(A \ifff B) \land D(A, r)) \impl (\exists r)D(B, r) & \textrm{(2), (3)}\\
(5) & D(A, r) \impl ((\forall \bar q)(A \ifff B) \impl (\exists r)D(B, r)) & \textrm{(4), \log{IPC}}\\
(6) & (\exists r)(D(A, r) \impl ((\forall \bar q)(A \ifff B) \impl (\exists r)D(B, r))) & \textrm{(5), R}\exists \\
(7) & (\forall \bar q)(A \ifff B) \impl ((\exists r)D(A, r) \impl (\exists r)D(B, r)) & \textrm{(6), \log{IPC}}
\end{eqnarray*}
The case of $C \equiv (\forall r)D(p, r)$ is handled similarly. \qed
\end{proof}

\begin{defn}
A formula $A$ of $\QG$ is in $\oper$-{\em normal form} if
it is quantifier-free and for all subformulas $\oper B$ of $A$, $B
\in \{\bot, \top\} \cup \Var$ or $B \equiv \oper B'$.
\end{defn}

\begin{proposition}\label{opernf}
Let $A$ be a quantifier-free formula of $\QG$. Then there exists a
formula $A'$ of $\QG$ in $\oper$-normal form such that
$\prov A \ifff A'$.
\end{proposition}

\begin{proof}
Follows from axiom (G1), Prop.~\ref{basicpr}(2) and (3)
using Prop.~\ref{equivsubst}(2).\qed
\end{proof}

\begin{proposition}
For every  $n \geq 0$,  $\prov \oper^n \top \ifff \top.$
\end{proposition}

\begin{proof}
$\oper^n \top \impl \top$ is already derivable intuitionistically.
For $\top \impl \oper^n \top$, use (G2), Prop.~\ref{basicpr}(1), and
induction on~$n$.\qed
\end{proof}

For propositional G\"odel logic, a normal form similar to the
disjunctive normal form of classical logic has been introduced in
\cite{baz96} (see also \cite{baaz-veith-98,baaz-veith-98b}).  This
so-called {\em chain normal form} is based on the fact that, in a
sense, the truth value of a formula only depends on the ordering of
the variables occurring in the formula induced by the valuation under
consideration.  The chain normal form can then be constructed by
enumerating all such orderings (using $\prec$ and $\ifff$ to encode
the ordering) in a way similar to how one constructs a disjunctive
normal form by enumerating all possible truth value assignments.  We
extend the notion of chain normal form and the results of
\cite{baaz-veith-98} in order to deal with the $\oper$ connective.
This is possible, since by Prop.~\ref{opernf} we can always push the
$\oper$ in front of atomic subformulas, so we only need to consider
orderings of subformulas of the form $\oper^j B$ with $B$ atomic. Let
$\Gamma$ be a finite subset of $\{ \oper^j p, \oper^j \bot : p \in
\Var, j \in \omega\} \cup \{\top\}$ and $\top, \bot \in \Gamma$.

\begin{defn}
A $\oper$-{\em chain} over $\Gamma$ is an expression  of the form
\[
(S_1 \star_1 S_2) 
\land  \cdots \land (S_{n-1} \star_{n-1} S_{n}) 
\]
such that $\Gamma = \{S_1, \ldots, S_n\}$, $S_1 \equiv \bot$, $S_n
\equiv \top$, and $\star_i \in \{\ifff, \prec\}$, for all $i=1, \dots
,n$.
\end{defn}  

Every $\oper$-chain $C$ uniquely determines a partition $\Pi_1^C$,
\dots, $\Pi_k^C$ of $\Gamma$ so that $\Pi_i^C = \{S_{j_i}, \ldots,
S_{j_{i+1}-1}\}$ where $j_1 = 1$, $j_{k+1} = n+1$, $j_{i} < j_{i+1}$,
$\star_{j_i} = \cdots = \star_{j_{i+1}-2} = {\ifff}$, and
$\star_{j_{i+1}-1} = {\prec}$.  Conversely, every such partition
determines a $\oper$-chain up to provable equivalences.  It is easily
seen that if $C$ is such a chain, then $\prov C \impl (S_i \ifff S_j)$
if $S_i, S_j \in \Pi_l^C$ for some $l$, and $\prov C \impl (S_i \prec
S_{i'})$ if $S_i \in \Pi_j^C$, $S_{i'} \in \Pi_{j'}^C$ and~$j <
j'$. Thus $C$ also uniquely corresponds to an ordering of $\Gamma$
which we denote $<_C$, defined by $S_{i} <_C S_{i'}$ iff $S_i \in
\Pi^C_j$, $S_{i'} \in \Pi_{j'}^C$ and $j < j'$. This order is total,
the $\Pi_i^C$ are maximal anti-chains, $\bot$ is minimal, and $\top$
is maximal.

Suppose now that $A$ is in $\oper$-normal form, and that $\Gamma$
contains all the subformulas of $A$ of the form $\oper^j p$ or $\oper
^j \bot$, as well as $\top$; that $C$ is an $\oper$-chain on $\Gamma$;
and that the valuation $v$ agrees with $<_C$, i.e., $S_i <_C S_j$ iff
$v(S_i) < v(S_j)$.  Using the same idea as in the proof of Lemma~3 in
\cite{baaz-veith-98}, one can find $A^C \in \Gamma$, the ``value'' of
$A$ under $C$, so that $v(A^C) = v(A)$, and the choice of $A^C$
depends only on $<_C$, not on $v$ itself.  Specifically, $A^C$ can be
constructed as follows: (1) If $A \in \Gamma$, then $A^C \equiv
A$. (2) If $A \equiv D \land E$, then $A^C \equiv D^C$ if $D^C <_C
E^C$ and $\equiv E^C$ otherwise. (3) If $A \equiv D \lor E$, then $A^C
\equiv D^C$ if $E^C <_C D^C$, and $\equiv E^C$ otherwise.  (4) If $A
\equiv D \impl E$, then $A^C \equiv E^C$ if $E^C <_C D^C$, and $\equiv
\top$ otherwise.  This ``evaluation'' of $A$ is provable in the sense
that $\prov C \impl (A \ifff A^C)$.  This follows easily using
the following theorems of~\LC:
\[\begin{array}{rcl@{\qquad}rcl}
(D \prec E) & \impl & (D \land E \ifff D) &
(E \prec D) & \impl & (D \land E \ifff E) \\
(D \ifff E) & \impl & (D \land E \ifff D) &
(D \prec E) & \impl & (D \lor E \ifff E) \\
(E \prec D) & \impl & (D \lor E \ifff D) &
(E \ifff D) & \impl & (D \lor E \ifff E) \\
(D \prec E) & \impl & (D \impl E \ifff \top) &
(E \prec D) & \impl & (D \impl E \ifff E) \\
(E \ifff D) & \impl & (D \impl E \ifff \top) 
\end{array}
\]
\begin{defn}
Let $A$ be a quantifier free formula in $\oper$-normal form,
$\Gamma_A$ be the set of all subformulas of $A$ of the form
$\oper^j p, \oper^k \bot, \top$, $\Gamma \supseteq \Gamma_A$, and
$C(\Gamma)$ the set of all possible $\oper$-chains over $\Gamma$. Then
\[
\bigvee_{C \in C(\Gamma)} C \land A^{C}
\]
is the $\oper$-{\em chain normal form} for $A$ over $\Gamma$. 
\end{defn}

\begin{theorem}\label{chainnf}
Let $A$ and $\Gamma$ be as above, and $A'$ be the $\oper$-chain normal
form for $A$ over $\Gamma$. Then $\prov A \ifff A'$.
\end{theorem}

\begin{proof} (See also Thm.~4 of \cite{baaz-veith-98}.)
First note that $\bigvee_{C \in C(\Gamma)} C$ is a tautology and provable
in~\log{LC}.  Since for each $C \in C(\Gamma)$ we have $\prov (C \land A^C)
\impl A$, the right-to-left implication $A' \impl A$ follows by case
distinction.

For the left-to-right implication, consider $A \impl (A \land \bigvee_{C \in
C(\Gamma)} C)$.  This is provable, since $\bigvee_{C \in
C(\Gamma)} C$ is provable. By distributivity of $\land$ over $\lor$,
we have $A \impl \bigvee_{C \in C(\Gamma)} (A \land C)$. We also have
$(A \land C) \impl (C \land A^C)$ for each $C \in C(\Gamma)$ from
$\prov C \impl (A \ifff A^C)$.  Together we get $A \impl \bigvee_{C \in
C(\Gamma)}(C \land A^C).$\qed
\end{proof}

We now strengthen the $\oper$-normal form result so that only
$\oper$-chains that are intuitively ``possible'' need to be
considered.  For this, we have to verify that we can exclude
chains~$C$ which result in orders which, e.g., have $\oper S <_C S$.

\begin{defn}
\label{minnorform}
A formula $A$ is in {\em minimal normal form} over $\Gamma$ if it is
of the form $\bigvee_{C \in {\cal C} \subseteq C(\Gamma)} C$, where each 
$C$ is a
$\oper$-chain over $\Gamma$, and so that the corresponding ordered
partition $\Pi_1^{C}, \ldots, \Pi_k^{C}$ satisfies
\begin{enumerate}
\item for no $i <j$ and $S \in \Gamma$ do we have $\oper^{r + s} S \in \Pi_i^C$
and $\oper^{r} S \in \Pi_j^{C}$ with $s > 0$;
\item for all $S \in \Gamma$, if $\oper^s S \in \Pi_i^C$ ($i < k$), then
$\oper^r S \notin \Pi_i^C$ if $r \neq s$; and
\item for no $j, j'$ and $S \in \Gamma$ do we have both $\oper^i S \in \Pi_j^C$
and $\oper^{i+1} S \in \Pi_{j'}^C$ with $j' > j + 1$. 
\end{enumerate}
\end{defn}

\begin{theorem}\label{minnf}
Let $A$ be in $\oper$-normal form.  There exists a formula $A^\nf$ in
minimal normal form such that $\prov A \ifff A^\nf$.
\end{theorem}

\begin{proof}
By Thm.~\ref{chainnf}, $\prov A \ifff A'$ where $A'$ is a $\oper$-chain
normal form over $\Gamma$.  Consider a disjunct of $A'$ of the form $C
\land A^C$, where $\Pi_1^C$, \dots, $\Pi_k^C$ is the ordered
partition of $\Gamma$ corresponding to $C$.  If $A^C \in \Pi_k^C$, then
$\prov (C \land A^C) \ifff C$, since $\prov A^C \ifff 
(A^C \ifff \top)$.  
Otherwise, $A^C \in \Pi_i^C$ with $i < k$. 
Then the sequence $\Pi_i^C$, \dots, $\Pi_k^C$ corresponds to a conjunction
\[
C' \equiv (A^C \star_{1} S'_1) \land \ldots \land (S'_{j-1} \star_j \top)
\]
where for at least one $l \le j$, $\star_j = \prec$, and $\prov C
\ifff C'' \land C'$, where $C''$ is the part of $C$ corresponding to
$\Pi_1^C$, \dots, $\Pi_{i-1}^C$.  
Since $\prov A^C \ifff (A^C \ifff \top)$, we have
\begin{equation}
\prov (C' \land A^C) \ifff (C' \land (\top \ifff  A^C)) \label{eqn1}
\end{equation}
As is easily seen, the right-hand side of~(\ref{eqn1}) is provably
equivalent to
\[
C'''\equiv (A^C \ifff S'_1) \land \ldots \land (S'_{j-1} \ifff \top)
\]
In sum, $\prov (C \land A^C) \ifff (C'' \land C''')$, and
$C'' \land C'''$ is a $\oper$-chain.

By induction on the number of disjuncts in $A'$ one shows that
there is $A''$ which is a disjunction of $\oper$-chains such that
$\prov A \ifff A''$. Now we have to prove that there exists a
disjunction of $\oper$-chains $A^\nf$ satisfying 1--3 of
Def.~\ref{minnorform} so that $\prov A'' \ifff A^\nf$.
 
Suppose that for some disjunct $C$ in $A''$ we have $\oper^{r+s} S \in
\Pi_i^C$ and $\oper^r S \in \Pi_j^C$ where $s > 0$ and $i < j$.  Then,
since $\prov (\oper^{r+s} A \prec \oper^r A) \ifff \oper^r A$ we have
$\prov C \ifff C'$ where $C'$ is the $\oper$-chain corresponding to
$\Pi_1^C, \dots, \Pi_{i - 1}^C, \Pi_i^C \cup \ldots \cup \Pi_k^C$.
 
Consider a disjunct $C$ of $A''$ where for some $i<k$, both
$\oper^r S\in \Pi_i^C$ and $\oper^s S \in \Pi_i^C$ where $r < s$.  Then
$\prov C \impl (\oper^s S \ifff \top)$.  To see this, recall that $\prov
\oper^r v \prec \oper^s S$ if $r < s$.  By definition of $\prec$, that
means that 
\begin{equation}
\prov ((\oper^s S \impl \oper^r S) \impl \oper^r S) \land
(\oper^r S\impl \oper^s S). \label{rsprec}
\end{equation}
Since $\prov C \impl (\oper^s S \ifff \oper^r S)$, we have $\prov C
\impl (\oper^s S \impl \oper^r S)$ which together with the left
conjunct of (\ref{rsprec}) gives $\prov C \impl \oper^r S$.  Thus, as
before, $C$ is provably equivalent to the $\oper$-chain corresponding to
$\Pi_1^C$, \dots, $\Pi_i^C \cup \ldots \cup \Pi_k^C$.  

Lastly, suppose that for a disjunct $C$ of $A''$ we have both $\oper^i
S\in \Pi_j^C$ and $\oper^{i+1} S \in \Pi_{j'}^C$ for some $j$, $j'$ such
that $j' > j+1 $. Then by axiom (G6) together with transitivity we get
$C \impl (\oper^{i+1} S \prec \oper^{i+1} S)$, and since $\prov (B
\prec B) \ifff B$ we have $\prov C \ifff C'$ where $C'$ is the
$\oper$-chain corresponding to $\Pi_1^C, \dots, \Pi_{j - 1}^C, \Pi_j^C
\cup \ldots \cup \Pi_{j'}^C \ldots \cup \Pi_k^C$.

By induction on the number of disjuncts in $A''$ we
obtain the desired $A^\nf$. \qed
\end{proof}
 

\section{Quantifier Elimination}

In this section we prove quantifier elimination for $\QG$.  As a
corollary of this result we show that the system $\QG$ is sound and
complete for $\G$ and that the latter is the intersection of all
finite-valued quantified propositional G\"odel
logics~$\mathbf{G}_k^\qp$.

\begin{proposition}\label{telprop}
\begin{enumerate}
\item $\prov (\forall p)A(p) \ifff  (A(\bot) \land (\forall p)A(\oper p))$ 
\item $\prov (\exists p)A(p) \ifff  (A(\bot)  \lor (\exists p)A(\oper p)).$  
\end{enumerate}
\end{proposition}

\begin{proof}
(1) The left-to-right implication follows easily from the two
    instances of (${\impl}{\forall}$)
\[
(\forall p)A(p) \impl A(\bot) \qquad\textrm{and}\qquad (\forall p)A(p)
\impl A(\oper p).
\]
For right-to-left, consider
\begin{eqnarray}
(q \ifff \bot) & \impl & (A(\bot) \land (\forall p)A(\oper p)) \impl A(q) \\ 
(q \ifff \oper p) & \impl & (A(\bot) \land (\forall p)A(\oper p)) \impl A(q) \label{eqn2}
\end{eqnarray}
which are derived easily from Prop.~\ref{equivsubst}(2) using $\log{IPC}^\qp$.
Use (R$\exists$) to introduce the existential quantifier in the antecedent of (\ref{eqn2}), and then (I12) to obtain
\begin{equation}
[(q \ifff \bot) \lor (\exists p)(q \ifff \oper p)] \impl  
(A(\bot) \land (\forall p)A(\oper p)) \impl A(q) \label{eqn3}
\end{equation}
The antecedent of (\ref{eqn3}) is an instance of (G5), and so 
\[
\prov (A(\bot) \land (\forall p)A(\oper p)) \impl A(q) 
\]
from which the right-to-left direction of~(1) follows by
(R$\forall$).  

(2) The argument is analogous to the derivation
of~(1).\qed
\end{proof}

\begin{defn}
For $\Gamma \subseteq \Var\cup\{\bot, \top\}$, let $\OP_\Gamma(A)$ be
the set of formulas inductively defined as follows:
\begin{eqnarray*}
\OP_\Gamma(A \ast B) &=& \OP_\Gamma(A) \cup \OP_\Gamma(B), \quad  {\rm where} 
\; \ast \in \{\vee, \wedge, \impl\}\\
\OP_\Gamma ((\Q p)A) & = & \OP_\Gamma(A), \quad  {\rm where} \; 
\Q \in\{\forall, \exists\}\\
\OP_\Gamma(\oper^k v)& = & 
                      \begin{cases}
                     \{\oper^k v\} &  {\rm if} \; v \in  \Gamma\\
                    \emptyset & {\rm otherwise}
                     \end{cases}
\end{eqnarray*}
Then $\exp_\Gamma(A) = \{k : \oper^k q \in \OP_\Gamma(A)\}$ 
\end{defn}

\begin{defn} The {\em
quantifier depth}~$\qd(A)$ of a formula is defined by:
\[\begin{array}{c@{\qquad}c}
\qd(p) = \qd(\bot)  =  0 & \qd((\forall p) B) = \qd((\exists p) B)  =  \qd(B) + 1 \\
\multicolumn{2}{c}{
\qd(B * C)  =  \max (\qd(B), \qd(C)) \textrm{ for } * \in \{\land, \lor, \impl\}}
\end{array}\]
\end{defn}

\begin{lemma}\label{mainlemma}
Let $A$ be a closed formula such that (a) every quantifier free
subformula of~$A$ is in $\oper$-normal form and (b) no two quantifier
occurrences bind the same variable. Let $\Delta = \{p_1, \ldots,
p_j\}$ be the set of variables belonging to the innermost quantifiers
in $A$, and $\Gamma = \Var(A) \setminus \Delta$.  Then there is a
formula $A^\sharp$ so that \begin{enumerate}
\item $\prov A\ifff A^\sharp$,
\item $\max \exp_\Delta(A^\sharp) \le \min \exp_\Gamma(A^\sharp)$,
\item $\max \exp_{\Var(A^\sharp)}(A^\sharp) \le 2\cdot \max\exp_{\Var(A)}(A)$, 
\item $\qd(A^\sharp) \le \qd(A)$.
\end{enumerate}
\end{lemma}

\begin{proof} Suppose $\Gamma = \{q_1, \ldots, q_l\}$. 
Let $A_0 = A$, $m = \max \exp_\Delta (A)$.  At stage $i$, pick the
non-innermost quantified subformula $(\forall q_i) B_i(q_i)$ or
$(\exists q_i) B_i(q_i)$ of $A_i$ corresponding to $q_i$ and replace
\begin{eqnarray*}
(\forall q_i) B_i(q_i) & \textrm{\ by\ } & B_i(\bot) \land \ldots \land
B_i(\oper^{{m}-1} \bot) \land (\forall p) B_i(\oper^{m} q_i) \\ (\exists q_i)
B_i(p) & \textrm{\ by\ } & B_i(\bot) \lor \ldots \lor
B_i(\oper^{{m}-1}\bot ) \lor (\exists q_i) B_i(\oper^{m} q_i)
\end{eqnarray*}
to obtain $A_{i+1}$.  The procedure terminates with $A_l = A^\sharp$.

At each stage $\prov A_i \ifff A_{i+1}$ follows by induction on $m$
from Prop.~\ref{telprop}. The lower bounds are obvious from the
construction of $A^\sharp$.\qed
\end{proof}

\begin{lemma}\label{sublemma}
Suppose $A(p)$ is in $\oper$-normal form and 
\[
\max \exp_{\{p\}} A \le \min\exp_{\Var(A)\setminus\{p\}}A.
\] 
There is a formula $A^\exists$,
with $\Var(A^\exists) \subseteq \Var(A) \setminus \{p\}$ so that
\[
\prov (\exists p) A \ifff A^\exists
\]
and $\max\exp_{\Var(A^\exists)\cup\{\bot\}} A^\exists \le
\max\exp_{\Var(A)\cup\{\bot\}}A + 1$.
\end{lemma}

\begin{proof}
Let $m = \max \exp_{\Var(A) \cup \{\bot\}} A$ be the maximal exponent
of a subformula $\oper^j S$ and let $\Gamma =
\{\oper^i S : S \in \Var \cup \{\bot\}, i \le m\}$.

Theorem~\ref{minnf} provides us with $A^\nf$ in minimal normal form
over $\Gamma$ so that $\prov (\exists p) A \ifff (\exists p) A^\nf$.
Since $\exists$ distributes over $\lor$, we only have to consider
formulas of the form $(\exists p) C$ where $C$ is a $\oper$-chain and
satisfies the conditions of Thm.~\ref{minnf}.
$C$ corresponds to an ordered partition $\Pi_1$, \ldots, $\Pi_k$ over
$\Gamma$.  We prove that $\prov (\exists p) C \ifff C'$ for some
quantifier-free $C'$ by induction on $k$.

If $k = 2$, then either $p \in \Pi_1$ or $p \in \Pi_k$.  In the first
case, $\prov (\exists p) C(p) \ifff C(\bot)$, in the second one, $\prov
(\exists p) C(p) \ifff C(\top)$.

Now suppose $k > 2$.  Three cases arise, according to how the 
equivalence classes containing $p$ are distributed.

(1) The partition corresponding to $C$ is of the form 
\[ 
\Pi_1, \ldots, \Pi_i, \{p\}, \{\oper p\}, \ldots, \{\oper^j
p\}\cup\Pi_k
\]
Then $C(p)$ is of the form 
\[ 
B \land \underbrace{(v \prec p) \land (p \prec \oper p) \land \ldots \land
(\oper^j p \ifff \top)}_{D(p)} \land \, E
\]
Since $D(\top)$ is provable, $\prov (\exists p) C \ifff B \land v \prec
\top \land E$.

(2) The partition corresponding to $C$ is of the form
\[ 
\Pi_1, \ldots, \Pi_i, \{p\}, \{\oper p\}, \ldots, \{\oper^j
p\}, \Pi_{i'}, \ldots, \Pi_k
\]
and $\oper^j p \notin \Pi_{i'}$. Then $C(p)$ is of the form 
\[ 
B \land \underbrace{(S \prec p) \land (p \prec \oper p) \land \ldots
\land (\oper^j p \prec S')}_{D(p)} \land E
\]
We first show that $\prov (\exists p) D(p) \ifff (\oper^{j+1} S \prec S')$.
For the right-to-left direction, observe that 
\[
\prov (\oper^{j+1} S \prec S') \impl [(S \prec \oper S) \land \ldots
\land (\oper^{j}S \prec \oper^{j+1}S) \land (\oper^{j+1} S \prec S'),
\]
from which the claim follows by~($\mathrm{R}{\exists}$).  The
left-to-right direction is proved by induction on $j$, using
axiom~(G6). In sum, we have
\[
\prov (\exists p) C(p) \ifff (B \land (\oper^{j+1} S \prec S') \land E)
\]

(3) The partition corresponding to $C$ is of the form
\[ 
\Pi_1, \ldots, \Pi_i, \{p\}, \{\oper p\}, \ldots, \{\oper^j
p\}\cup\Pi, \Pi_{i'}, \ldots, \Pi_k
\]
with $S \in \Pi$, $S \neq \oper^j p$.  Because of the condition on $\max
\exp_{\{p\}} A$ we can assume that $S \equiv \oper^n q$ with $n \ge
j$.

We proceed by induction on $j$.  If $j = 0$, then we have a conjunct
$p \ifff S$, and $(\exists p) C \equiv C(S)$.  Otherwise, we have a
conjunct $\oper^j p \ifff \oper^n q$ with $n \ge j$.  Using (G3), this
conjunct is provably equivalent to $(\oper^{j-1} p \ifff \oper^{n-1} q)
\lor (\oper^j p \land \oper^n q)$.  Hence, $C$ is equivalent to the
disjunction of two $\oper$-chains corresponding to
\begin{eqnarray*}
\Pi_1, \ldots, \Pi_i, \{p\}, \{\oper p\}, \ldots, & & \{\oper^{j-1}
p, \oper^{n-1} q\}, \Pi, \Pi_{i'}, \ldots, \Pi_k\\
\Pi_1, \ldots, \Pi_i, \{p\}, \{\oper p\}, \ldots, & & \{\oper^j
p\}\cup \Pi \cup \Pi_{i'} \cup \ldots\cup \Pi_k
\end{eqnarray*}
For the first $\oper$-chain, the maximum exponent of $p$ is smaller
and hence the induction hypothesis of the present subcase applies.
The second $\oper$-chain is shorter overall, and hence the induction
hypothesis based on number of equivalence classes applies.\qed
\end{proof}

\begin{lemma}
Let $A(p)$ be in $\oper$-normal form, and so that
\[
\max \exp_{\{p\}} A \le \min\exp_{\Var(A)\setminus\{p\}}A.\] 
There is a
formula $A^\forall$, with $\Var(A^\forall) \subseteq \Var(A) \setminus
\{p\}$ so that
\[
\prov (\forall p) A \ifff A^\forall
\]
and $\max\exp_{\Var(A^\forall) \cup \{\bot\}} A^\forall \le
\max\exp_{\Var(A)\cup\{\bot\}}A + 1$.
\end{lemma}

\begin{proof}
Let $A^\nf$ be the minimal normal form of $A$.  It is provably
equivalent to the formula obtained from $A^\nf$ by replacing each
element of a chain $S \prec S'$ by $\oper S \impl S'$.  By
distributivity then, $A \ifff A'$ where $A'$ is a conjunction of
disjunctions of implications of the form $\oper^i S \impl
\oper^j S'$.  Any such disjunct of the form $\oper^i p \impl
\oper^j p$ is provably equivalent to $\top$ if $i \le j$ (in which
case the entire disjunction can be deleted), or to $\top \impl \oper^j p$
if $i > j$.  The part of a disjunction in $A'$ containing $p$ thus 
can be assumed to be of the form
\[
\bigvee_i (D_i \impl \oper^{n_i} p) \lor \bigvee_j (\oper^{m_j} p \impl E_j)
\]
where $p \notin D_i, E_i$.  This, in turn, is equivalent to a
conjunction of disjunctions of the form
\[
\bigvee_i (D \impl \oper^{n_i} p) \lor \bigvee_j (\oper^{m_j} p \impl E)
\]
This can again be simplified by taking $n = \max \{n_i\}$ and $m =
\min \{m_j\}$, since $\prov (A \impl B) \lor (A \impl C) \ifff (A \impl C)$
if $\prov B \impl C$.

Since $\prov (\forall p) (A \land B) \ifff (\forall p) A \land
(\forall p) B$ and $\prov (\forall p)(A(p) \lor B) \ifff (\forall p)
A(p) \lor B$ if $p \notin B$, it suffices to show that a formula of
the form
\[
F \equiv (\forall p)(D \impl \oper^n p) \lor (\oper^m p \impl E))
\]
is equivalent to a quantifier free formula.  We distinguish three cases:

(1) $E \equiv \oper^k \top$, $k \ge 0$. Then $\prov (\oper^m p \impl
    E)$ and hence $\prov F \ifff \top$.

(2) $E \equiv \oper^k \bot$, $k < m$. Then $\prov (\oper^m p \impl
    E) \ifff E$, and hence $\prov F \ifff (A \impl \oper^n \bot) \lor
    E$.

(3) Since $\max \exp_{\{p\}} A \le \min\exp_{\Var(A)\setminus\{p\}}A$
    by assumption, this leaves only the case $E \equiv \oper^m S$.
    Then $\prov F \ifff (A \impl \oper^{n+1} S) \lor \oper^m S$.  The
    left-to-right implication is obvious by (${\impl}{\forall}$),
    instantiating $p$ by~$\oper S$.  For the right-to-left implication
    two cases arise:

(a) $n \le m$.  By (G4), we have $\prov (A \impl \oper^{n+1} S) \impl
[(A \impl \oper^{n} p) \lor (\oper^n p \impl \oper^n S)]$.
Furthermore,  $\prov (\oper^n p \impl \oper^n S)
\impl (\oper^m p \impl \oper^m S)$. In sum, we have
\[
[(A \impl \oper^{n+1} S) \lor \oper^m S] \impl [(A \impl \oper^n p)
\lor (\oper^m p \impl \oper^m S) \lor \oper^m S]
\]
Since $\prov \oper^m S \impl (\oper^m p \lor \oper^m S)$, we have 
$\prov [(A \impl \oper^{n+1} S) \lor \oper^m S] \impl F$.

(b) $n > m$. By (G2), $\prov \oper^m S \impl \oper^{n+1} S$, and so
$\prov [(A \impl \oper^{n+1} S) \lor \oper^{m} S] \impl (A \impl
\oper^{n+1} S]$. Using induction and (G4), it is easy to show that
\[
\prov (A \impl \oper^{n+1} S) \impl [\underbrace{(A \impl \oper^n p)
\lor \bigvee_{i=m}^{n-1}(\oper^{i+1} p \impl \oper^i p)}_D \lor
(\oper^m p \impl \oper^m S].
\]
Each of the disjuncts $\oper^{i+1} p \impl \oper^i p$ implies
$\oper^i p$, which in turn implies $A \impl \oper^n p$, so $\prov D
\impl (A \impl \oper^n p)$. In sum, we have again $\prov
[(A \impl \oper^{n+1} S) \lor \oper^m S] \impl F$.

The bound on $\max \exp_{\Var(A^\forall) \cup \{\bot\}} A$ follows by
inspection.\qed
\end{proof}

\begin{theorem}
For every closed formula $A$ of $\QG$ there exists a 
variable-free formula $A^\qf$ such that $\prov A \ifff A^\qf$, and
$\max \exp_{\{\bot\}} A^\qf\le 2^{\qd(A) + l}$ where $l = \max
\exp_{\Var(A) \cup \{\bot\}}$.
\end{theorem}

\begin{proof}
We may assume, renaming variables if necessary, that each variable in
$A$ is bound by only one quantifier occurrence.  By induction on
$\qd(A)$.  If $\qd(A) = 0$, there is nothing to prove.  If $\qd(A) >
0$, let $A^\sharp$ be as in Lemma~\ref{mainlemma}.  Replace each
innermost quantified formula $(\exists p) B$, $(\forall p) B$ by
$B^\exists$ or $B^\forall$, respectively.  The resulting formula $A'$
satisfies $\qd(A') \le \qd(A) - 1$ and $\max \exp_{\Var(A)
\cup\{\bot\}} A' \le 2\max \exp_{\Var(A) \cup \{\bot\}} A+1$.\qed
\end{proof}

\begin{proposition}
Let $A$ be variable-free, and in $\oper$-normal
form. Then either $\prov A \ifff \top$ or $\prov A \ifff \oper^k
(\bot)$ where $k \le \max \exp_{\{\bot\}} A = n$.
\end{proposition}

\begin{proof}
Consider the minimal normal form $A^\nf$ of $A$ over $\{\oper^k(\bot)
: k \le n\}$.  Each chain in $A^\nf$ is of one of two forms
\begin{eqnarray*}
C & = & (\bot \prec \oper(\bot)) \land (\oper(\bot) \prec
\oper\oper(\bot)) \land \ldots \land (\oper^{n-1} \bot \prec
\oper^n(\bot)) \\[-2ex] 
C_m & = & (\bot \prec \oper(\bot)) \land
(\oper(\bot) \prec \oper\oper(\bot)) \land \ldots \land (\oper^{m-1}
\bot \prec \oper^m(\bot)) \land \bigwedge_{k=m}^n \oper^k(\bot)
\end{eqnarray*}
$C$ is provable, so $\prov C \ifff \top$, and $\prov C_m \ifff
\oper^m(\bot)$.  So if $A^\nf$ contains $C$, then $\prov A \ifff
\top$, otherwise $\prov A \ifff \oper^k(\bot)$, where $k$ is the
maximum of $C_i$ occurring in $A^\nf$.\qed
\end{proof}

\begin{corollary}
Let $A$ be closed and not containing $\oper$.  Then either $\prov A$
or $\prov A \ifff \oper^k(\bot)$, where $k \le 2^{\qd(A)}$.
\end{corollary}

\begin{corollary}\label{soundcompl}
The calculus $\QG$ is complete for $\G$. 
\end{corollary}

\begin{proof}
If $\QG \not\vdash A$, then $\QG \vdash A \ifff \oper^k\bot$ for some
$k$.  Since $\G \nmodels \oper^k\bot$ for all~$k$, $\G \nmodels A$.
\end{proof}

\begin{theorem}
$\G$ is the intersection of all finite-valued quantified
propositional G\"odel logics.
\end{theorem}

\begin{proof}
\QG{} is sound for each finite-valued G\"odel logic, so $\G \subseteq
{\bf G}_k^\qp$ for each $k$.  Conversely, if $\G \not\models A$, then
$\prov A \ifff \oper^k(\bot)$ for some $k$.  Since \QG{} is sound for
${\bf G}_{k+2}$, we have ${\bf G}_{k+2} \not\models A$ as obviously
${\bf G}_{k+2} \nmodels \oper^k\bot$.
\end{proof}


\end{document}